\documentclass{article}
\usepackage[utf8]{inputenc}
\usepackage[total={4.9in, 7.7in}]{geometry}
\usepackage{amsmath}
\usepackage{mathtools}
\usepackage{amssymb}
\usepackage{amsthm}
\usepackage{amsfonts}
\usepackage{mathrsfs}
\usepackage{graphicx}
\usepackage{float}
\usepackage{caption}
\usepackage{fancyhdr}

\newtheorem{theorem}{Theorem}
\newtheorem{lemma}{Lemma}
\newtheorem{proposition}{Proposition}
\newtheorem{corollary}{Corollary}

\fancypagestyle{titlefooter}{
  \fancyhf{}
  
  \fancyfoot[C]{This research was supported by NSF grant no. 1950563.}}

\title{\tb{Graphs with prescribed radius, diameter, and center}\vspace{-4ex}}
\author{}
\date{}

\newcommand{\tb}[1]{\textbf{#1}}

\begin{document}
\maketitle
\begin{center} \large
\textbf{\scalebox{1.0909}{Kelly Guest}} \vspace{2mm} \\ 
Tuskegee University, kbraganguest@gmail.com \vspace{4mm}\\

\textbf{\scalebox{1.0909}{Andrew Johnson}} \vspace{2mm} \\ 
Kennesaw State University, ajohn301@students.kennesaw.edu \vspace{4mm}\\

\textbf{\scalebox{1.0909}{Peter Johnson}} \vspace{2mm} \\ 
Auburn University, johnspd@auburn.edu \vspace{4mm}\\

\textbf{\scalebox{1.0909}{William Jones}} \vspace{2mm} \\ 
State University of New York at Binghamton, wjones4@binghamton.edu \vspace{4mm}\\

\textbf{\scalebox{1.0909}{Yuki Takahashi}} \vspace{2mm}\\
Grinnell College, takahash@grinnell.edu \vspace{4mm}

\textbf{\scalebox{1.0909}{Zhichun (Joy) Zhang}} \vspace{2mm}\\
Swarthmore College, zzhang3@swarthmore.edu \vspace{4mm}

\end{center}

\thispagestyle{titlefooter}

\bigskip

\begin{abstract}
    Among other things, it is shown that for every pair of positive integers $r$, $d$, satisfying $1<r<d\leq 2r$, and every finite simple graph $H,$ there is a connected graph $G$ with diameter $d$, radius $r$, and center $H.$
\end{abstract}

\bigskip



\section{Introduction}

All graphs referred to will be finite and simple. The vertex and edge sets of a graph $G$ will be denoted $V(G)$ and $E(G)$, respectively. If $G$ is connected and $u,v \in V(G)$, $dist_G(u,v)$ is the length of a shortest walk in $G$ from one of $u,v$ to the other; a geodesic under the shortest-walk metric. As every shortest walk is a path, $dist_G(u,v)$ may also be formulated as the length of a shortest path in $G$ with end-vertices $u$ and $v$.

If $G$ is connected and $v \in V(G)$, the \textit{eccentricity} of $v$ in $G$, denoted $\varepsilon_G(v)$, is:
$$\varepsilon_G(v) = \max_{u \in V(G)}\{dist_G(u,v)\}.$$
The \textit{radius} of a connected graph $G$ is:
$$rad(G) = \min_{u\in V(G)}\{\varepsilon_G(u)\},$$
and its \textit{diameter} is:
$$diam(G) = \max_{u \in V(G)}\{\varepsilon_G(u)\}.$$
Equivalently, \[diam(G) = \max_{u,v \in V(G)}\{dist_G(u,v)\}.\]

It is easy to see that $rad(G) \leq diam(G) \leq 2rad(G)$. It is a standard exercise in a first course in graph theory to show that for any positive integers satisfying $r \leq d \leq 2r$, there is a connected graph $G$ such that $rad(G) = r$ and $diam(G) = d$. (A more challenging, but still elementary, exercise would be to determine, for pairs $r,d$ constrained as above, the values of $n$ such that there exists a connected graph $G$ with $rad(G) = r$, $diam(G) = d$, and $|V(G)| = n$.)

A vertex $v \in V(G)$ is a \textit{central vertex} in $G$ if and only if $\varepsilon_G(v) = rad(G)$. The \textit{center} of $G$, denoted $C(G)$, is the subgraph of $G$ induced by the set of centers of $G$. (Therefore, that set is $V(C(G))$.)

The question broached in [1] is: which graphs can be installed as the center of another graph? That is, given a graph $H$, can you find a connected graph $G$ such that $C(G) \cong H$?

As reported in [1], this question in full generality was killed at its birth as a question meriting research by a brilliant observation of S. Hedetniemi (Steve? Sandra?), encapsulated in Figure 1.

\begin{figure}[H]
\centering
\captionsetup{justification=centering,margin=2cm}
\includegraphics[scale = 0.15]{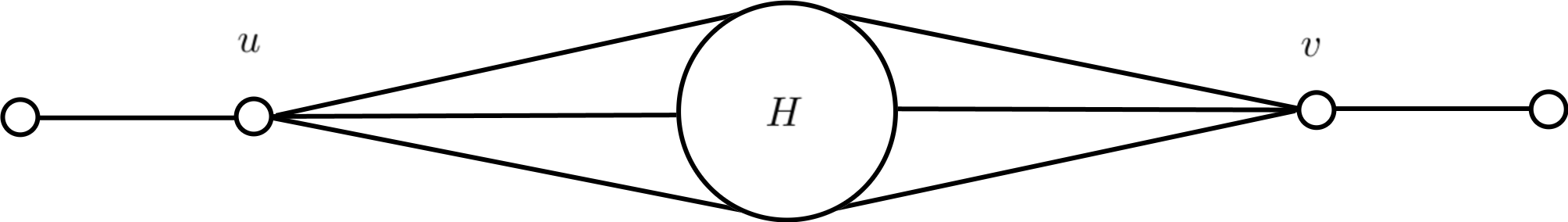}
\caption{A connected graph $G$ with an arbitrary graph $H$ as its center. Each vertex of $H$ is adjacent to both $u$ and $v$, in $G$.}
\end{figure}

The authors of [1] resurrect the problem by asking: for a distinguished family $\mathcal{F}$ of connected graphs, which graphs $H$ can be the center of a graph $G \in \mathcal{F}$? And, for such $H$ and $\mathcal{F}$, how small can $|V(G)| - |V(H)|$ be, if $G\in \mathcal{F}$?
These questions have borne fruit, but we are going in a different direction.

The graph $G$ in Figure 1 has diameter 4 and radius 2. The set of central vertices of $G$ is precisely $V(H)$, regardless of what $H$ is. If the paths leading away from $H$ from $u$ and $v$ are each lengthened to have length $t > 1$, the result is a graph with center $H$, radius $t+1$, and diameter $2t+2$.

Our aim here is to answer the question: for which positive integers $d,r$, satisfying $r \leq d \leq 2r$, and graphs $H$, does there exist a connected graph $G$ such that $rad(G) = r$, $diam(G) = d$, and $C(G) \simeq H$? The extension of the observation of Hedetniemi just above shows that there is such a $G$ for every $H$, $r>1$, and $d=2r$. Our main result, in Section 3, is that there is such a $G$ for every $H$, $r>1$, and $r<d\leq 2r$. In the next section we deal with extremes, and alternative solutions to that in Section 3, in some cases.

\section{Extremes and alternative solutions}

\subsection{$\boldsymbol{r = d}$}
If $rad(G) = diam(G)$, then $G$ is its own center. Therefore, $H = C(G)$ and $rad(G) = diam(G)$ if and only $H \simeq G$ and $rad(H) = diam(H)$.

\subsection{$\boldsymbol{r =1, d=2}$}
If $rad(G) = 1$, then each central vertex of $G$ is adjacent to every other vertex of $G$. Therefore, if $H \cong C(G)$ then $H$ must be a complete graph, and each vertex of $H$ must be adjacent to each vertex of $V(G) \setminus V(H)$. Furthermore, since all central vertices of $G$ are in $V(H)$, it must be that every $v \in V(G) \setminus V(H)$ has a non-neighbor in $G$ in $V(G) \setminus V(H)$.

Let ``$\vee$" stand for the \textit{join} of two graphs: $X \vee Y$ is formed by taking disjoint copies of $X$ and $Y$ and then adding in every edge $xy$, $x \in V(X), y \in V(Y)$. By the paragraph above, when $r=1$, $d=2$, the only $H$ for which a solution $G$ can exist are $H = K_t, t>0$, and the only possible solutions are $K_t \vee Y$ in which $Y$ is a graph with $|V(Y)| > 1$ and for each $y \in V(Y)$, the degree deg$(y)$ of $y \in V(Y)$ satisfies deg$_Y(y) < |V(Y)| - 1$.

Every such $G = K_t \vee Y$ satisfies $rad(G) = 1$, $diam(G) = 2$, and $C(G) = K_t$, so we have completely characterized the values of $H (H = K_t)$ for which our problem with $r=1, d=2$ has a solution, and all possible solutions ($G = K_t \vee Y$, as above).

\subsection{A standard method}

\begin{proposition} Suppose that $X$ is a connected graph with $|V(X)| >1$, $rad(X) > 1$, and $V(C(X)) = \{h\}$; i.e., there is a single central vertex in $X$. For an arbitrary graph $H$, if $G$ is formed by replacing $h$ by $H$, with every vertex of $H$ adjacent in $G$ to every vertex in $X$ to which $h$ is adjacent, then $rad(G) = rad(X)$, $diam(G) = diam(X)$, and $C(G) \cong H$.\end{proposition}

The proof is straightforward. Note that the assumption that $rad(X) = \varepsilon_X(h) \geq 2$ plays a role in the proof that $H \cong C(G)$. 

For instance, the graph in Figure 1 is obtained from $X = P_5$, the path on 5 vertices, by the device of Proposition 1. The generalization to the solution of our problem for all $H$ when $d = 2r \geq 4$ uses the device of Prop. 1 with $X = P_{2r+1}$.

In Figure 2 we have a graph $X$ with a single central vertex $h$ such that $rad(X) = r$, $diam(G) = 2r-1$, for arbitrary $r \geq 2$. By Proposition 1, this shows that every graph $H$ can be the center of a graph $G$ of radius $r$ and diameter $2r-1$, for every $r \geq 2$.

\begin{figure}[H]
\centering
\captionsetup{justification=centering,margin=2cm}
\includegraphics[scale = 0.185]{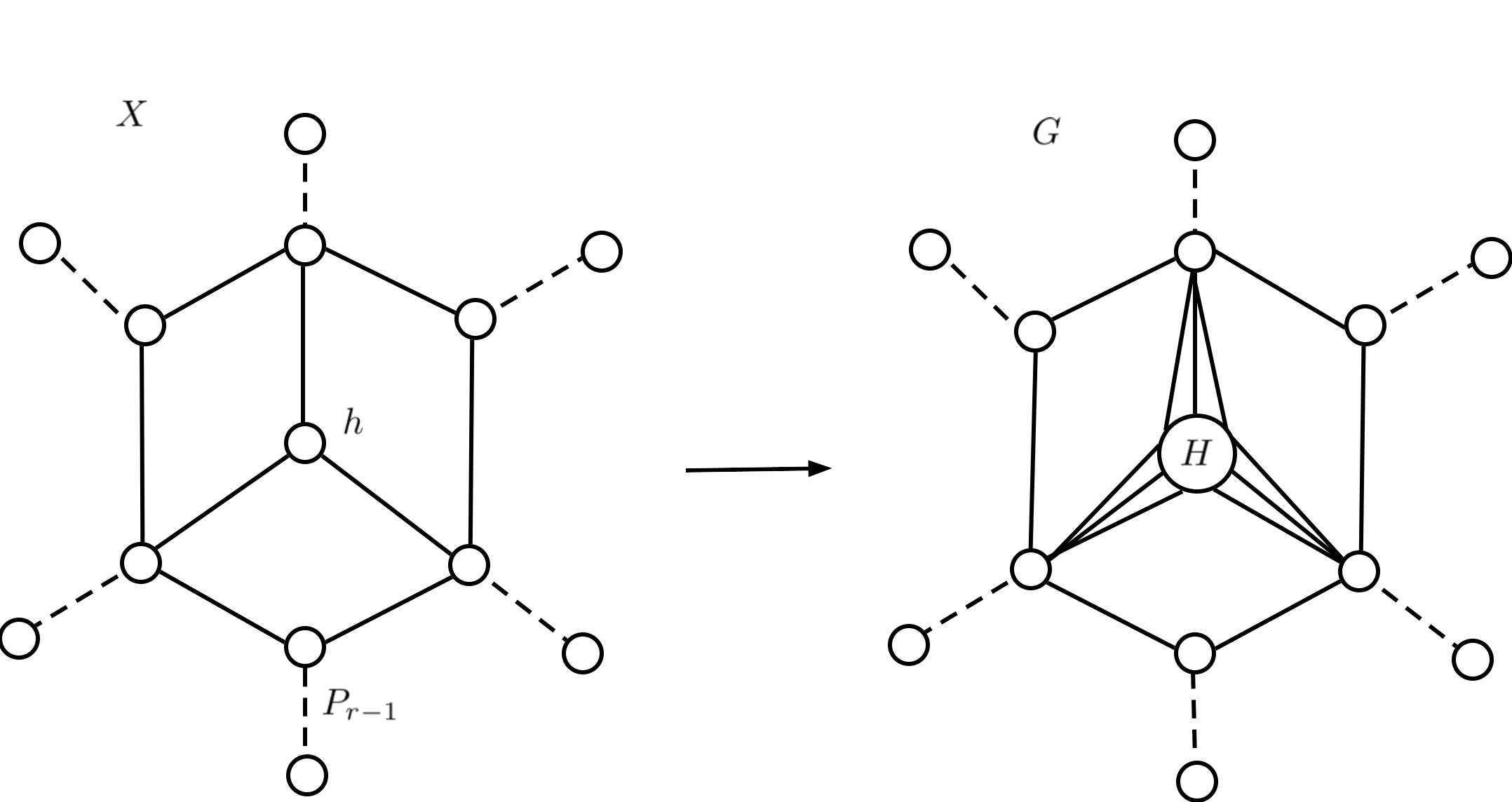}
\caption{A graph $X$ with radius $r\geq 2$, diameter $2r-1$, and a single central vertex $h$; and a graph $G$ with $rad(G)=r$, $diam(G)=2r-1$, and $C(G)\simeq H$. The paths hanging off the vertices of $C_6$ are all $P_{r-1}$, paths of length $r-2$. In the case $r=2$, they are not there, and $|V(X)|=7$.}
\end{figure}

For those who enjoy variety, we can vary $X$ to the graph $Y$ shown in Figure 3, which gives another solution to our problem when $d=2r$ and $H$ arbitrary. 

\begin{figure}[H]
\centering
\captionsetup{justification = centering, margin = 2cm}
\includegraphics[scale = 0.185]{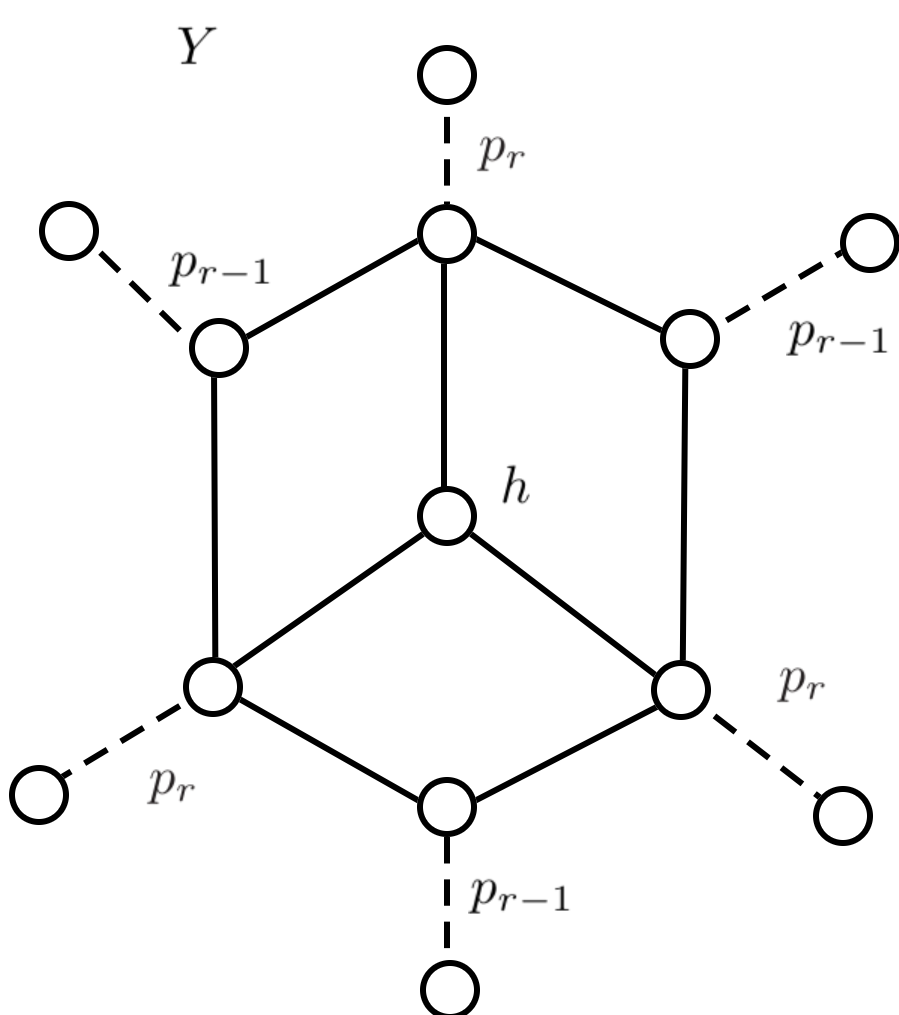}
\caption{A graph with a single central vertex, radius $r\geq 2$, and diameter $2r$.}
\end{figure}

If you have been paying attention, you might exclaim: why do we need this? Hedetniemi's construction already gives us solutions of our problem in the case $d = 2r \geq 4$. Yes, bur Figure 3 gives a \textit{different} solution, and different solutions of our problem contribute to the solution of a problem that towers over ours: given positive integers $r$ and $d$ satisfying $1 < r < d \leq 2r$, and a graph $H$, find all possible graphs $G$ satisfying $rad(G) = r$, $diam(G) = d$, and $C(G) \cong H$. In view of Proposition 1, in pursuit of this larger problem, it is appropriate to pose the following: given $d$ and $r$ as above, find all graphs $X$ such that $rad(X) = r$, $diam(X) = d$, and $C(X) = K_1$.

Moreover, the alternative solutions to the $d=2r$ case provide a related problem: what properties characterize those graphs with $d=2r$ and center $K_1$? The majority of graphs constructed with center $K_1$ in fact had $d=2r$, and the solution to this problem will considerably narrow down the larger problem.

In Figure 4, we have, for $r \geq 2$, a graph of radius $r$ and diameter $r+1$, and a graph of radius $r$ and diameter $r+\lceil\frac{r}{3}\rceil$, both with a single central vertex.

\begin{figure}[H]
\centering
\captionsetup{justification=centering,margin=2cm}
\includegraphics[scale = 0.15]{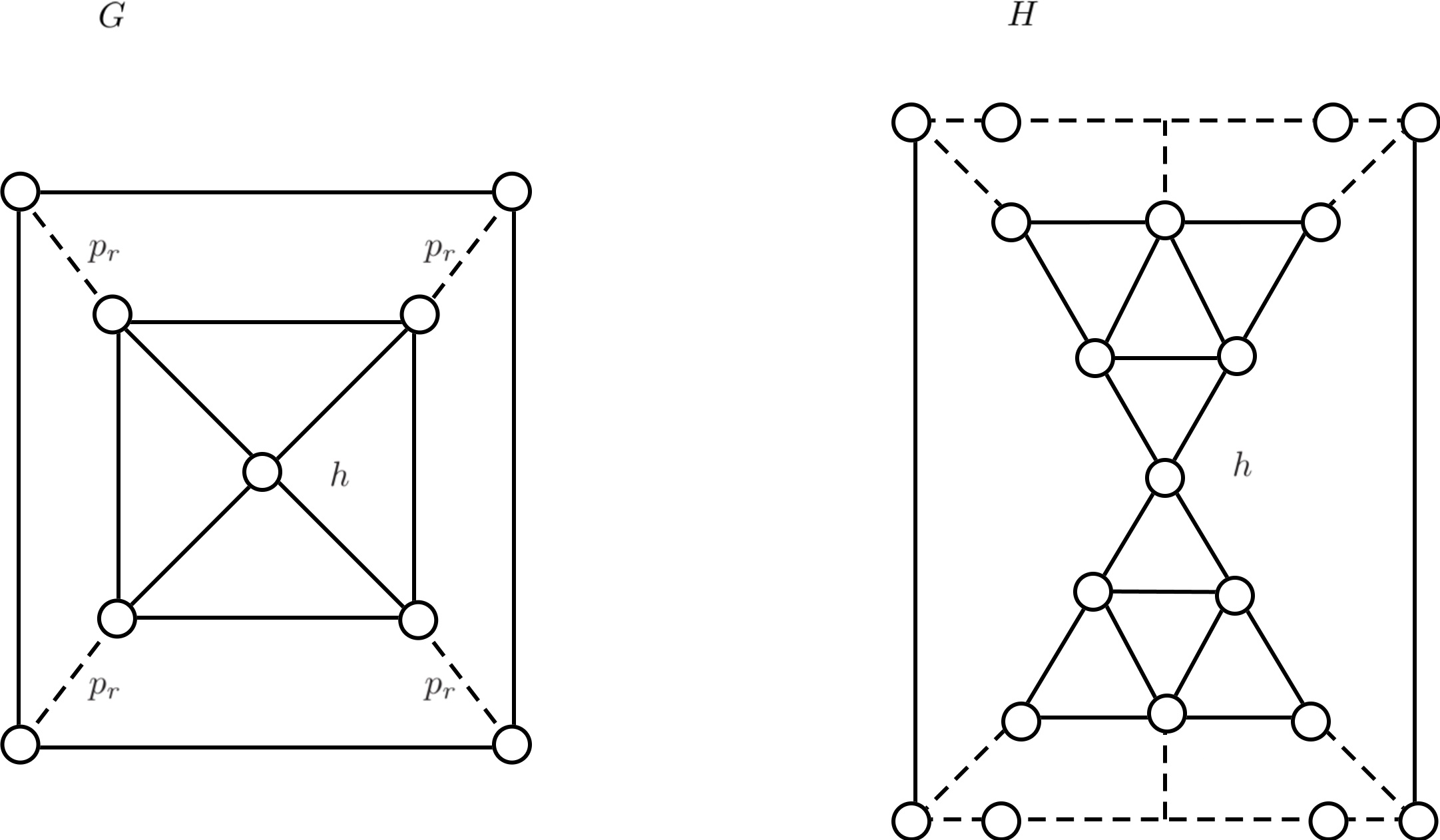}
\caption{A graph $G$ with radius $r$ and diameter $r+1$, and a graph $H$ with radius $r$ and diameter $r+\lceil\frac{r}{3}\rceil$.
The "top" and the "bottom" of the drawing of $H$ are $P_{r+1}$'s.}
\end{figure}

\subsection{A non-standard strategy in special cases}

The strategy referred to, applicable only when $H$ is connected is: attach pairwise vertex-disjoint paths to the vertices of $H$. This trick appears to be of use only in a special class of cases.

\begin{proposition}
Suppose that $H$ is connected with $rad(H) = diam(H) = z$. Suppose that $G$ is formed by attaching vertex-disjoint paths $P_t$ to the vertices of $H$, with each vertex of $H$ being an end of its attached path (when $t=1$, nothing is attached, and $G = H$). Then $rad(G) = z+t-1$, $diam(G) = 2(t-1) +z$, and $C(G) \cong H$.
\end{proposition}

The proof is straightforward.

\begin{corollary}
If $H$ is as in Proposition 2 then for all integers $r\geq z$ and $d = 2r-z$ there is a graph $G$, obtained as in Prop. 2 with $t=r-z+1$ such that $rad(G) = r$, $diam(G)=d$, and $C(G) = H$.
\end{corollary}

\section{The main result}

\begin{lemma}
Let $X$ be the graph depicted in Figure 5. Suppose that $n \geq 0$ and $r \geq \max\{2,n+1\}$. Then $h$ is the unique central vertex of $X$, $rad(X) = r$, and $diam(X) = r+n+1$.
\end{lemma}

\begin{figure}[h]
\centering
\captionsetup{justification=centering,margin=2cm}
\includegraphics[scale = 0.11]{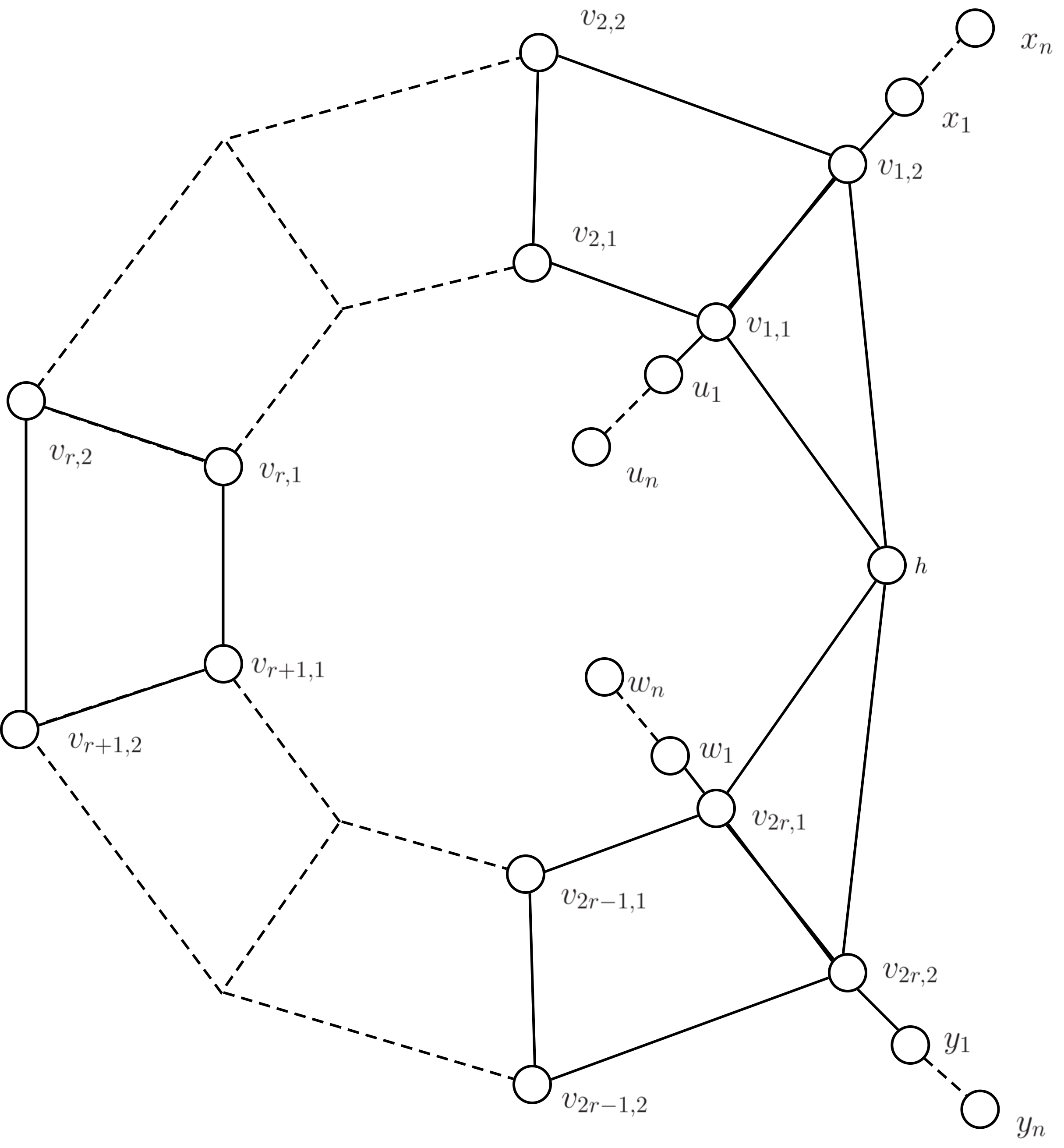}
\caption{A graph $X$ with a single central vertex $h$, radius $r$ and diameter $r+n+1$, provided $r\geq n+1$.}
\end{figure}

\begin{proof}
Clearly $\varepsilon_X(h) = \max\{r,n+1\}=r$. Checking shows that every other vertex of $X$ has eccentricity $>r$ in $X$. For instance, \[\varepsilon_X(v_{1,1}) = \max\{dist(v_{1,1},w_n),dist(v_{1,1},v_{r+1,2})\} = \max\{n+2,r+1\} = r+1.\]
Finally, it is easy to see that the vertices $v_{i,j}$, $i \in \{r,r+1\}$, $j\in \{1,2\}$, have the greatest eccentricity; for instance, $\varepsilon_X(v_{r,1}) = dist(v_{r,1},y_n) = r+n+1 = diam(X)$.
\end{proof}

\begin{theorem}
For all integers $r\geq 2$ and $d$ satisfying $r < d \leq 2r$ and every graph $H$ there is a graph $G$ such that $rad(G) = r$, $diam(G) = d$, and $C(G) \cong H$. Furthermore, $G$ is obtainable from some graph by the method of Proposition 2.
\end{theorem}

\section*{Acknowledgement}
We greatly thank Timothy Eller for drawing the graphs for us.


\begin{thebibliography}{9}
\bibitem{latexcompanion} 
Fred Buckley, Zevi Miller, and Peter J. Slater. On graphs containing a given graph as center.
\textit{Journal of Graph Theory}
5(1981), 427-434.

\end{thebibliography}
\end{document}